\documentclass{amsart}
\usepackage{amssymb,amsmath}

\newtheorem{theorem}{Theorem}
\newcommand{\bt}{\begin{theorem}}
\newcommand{\et}{\end{theorem}}

\newtheorem{lemma}{Lemma}
\newcommand{\bl}{\begin{lemma}}
\newcommand{\el}{\end{lemma}}

\newtheorem{example}{Example}
\newcommand{\bex}{\begin{example}}
\newcommand{\eex}{\end{example}}

\newtheorem{problem}{Problem} 
\newcommand{\bp}{\begin{problem}}
\newcommand{\ep}{\end{problem}}

\newcommand{\beq}{\begin{equation}}
\newcommand{\eeq}{\end{equation}}
\newcommand{\benum}{\begin{enumerate}}
\newcommand{\eenum}{\end{enumerate}}
\newcommand{\ba}{\begin{array}}
\newcommand{\ea}{\end{array}}
\newcommand{\Z}{\ensuremath{\mathbf Z}}
\newcommand{\N}{\ensuremath{\mathbf N}}
\newcommand{\Q}{\ensuremath{\mathbf Q}}

\newcommand{\C}{\ensuremath{\mathbf C}}

\newcommand{\mcp}{\ensuremath{\mathcal P}}

\newcommand{\F}{\ensuremath{\mathbf F }}
\newcommand{\Fp}{\ensuremath{ {\mathbf F}_p }}

\DeclareMathOperator{\card}{card}

\DeclareMathOperator{\Gal}{Gal}

\begin{document}
\date{\today}
\title{Binary linear forms over finite sets of integers}
\author[M. B. Nathanson]{Melvyn B. Nathanson}
\address{Department of Mathematics\\Lehman College (CUNY)\\Bronx, New York 10468}
\email{melvyn.nathanson@lehman.cuny.edu}
\thanks{The work of M.B.N. was supported in part by grants from the NSA Mathematical Sciences Program
and the PSC-CUNY Research Award Program.}

\author[K. O'Bryant]{Kevin O'Bryant}
\address{Department of Mathematics\\ 
College of Staten Island (CUNY)\\Staten Island, New York 10314}
\email{kevin@member.ams.org}

\author[B. Orosz]{Brooke Orosz}
\address{Department of Mathematics\\ CUNY Graduate Center\\
New York, New York 10036}
\email{borosz@gc.cuny.edu}

\author[I. Ruzsa]{Imre Ruzsa}
\address{R\' enyi Institute of Mathematics\\
Hungarian Academy of Sciences\\
Budapest, Hungary}
\email{ruzsa@renyi.hu}

\author[M. Silva]{Manuel Silva}\address{Department of Mathematics\\New University of Lisbon\\Lisbon, Portugal}\email{mnasilva@gmail.com}


\begin{abstract}
Let $A$ be a finite set of integers.  For a polynomial $f(x_1,\ldots,x_n)$ with integer coefficients, let $f(A) = \{f(a_1,\ldots,a_n) : a_1,\ldots,a_n \in A\}.$  In this paper it is proved that for every pair  of normalized binary linear forms $f(x,y)=u_1x+v_1y$ and $g(x,y)=u_2x+v_2y$ with integral coefficients, there exist arbitrarily large finite sets of integers $A$ and $B$ such that $|f(A)| > |g(A)|$ and  $|f(B)| < |g(B)|$.
\end{abstract}

\date{\today}

\maketitle

\section{Polynomials over finite sets of integers}
Let $f = f(x_1,\ldots,x_n)$  be a polynomial with integral coefficients.
Let $A$ be a nonempty finite set of integers or of congruence classes modulo $m$.
We denote by $f(A)$ the image of the function $f$ with domain $A$, that is,
\[
f(A) = \{f(a_1,\ldots,a_n) : a_i \in A \text{ for } i = 1,\ldots, n\}.
\]
If $f$ is a polynomial in $n$ variables, then $|f(A)| \leq |A|^n.$

The classical examples are the polynomials $x_1+x_2$ and $x_1-x_2.$  The sumset $A+A$ is the set $s(A)$ for the polynomial $s(x_1,x_2) = x_1+x_2,$ and the difference set $A-A$ is the set $d(A)$ for the polynomial $d(x_1,x_2) = x_1-x_2.$  
For any arithmetic progression $A$ or, more generally, any symmetric set $A$ of integers, we have $|d(A)| = |s(A)|,$ but for ``most''  sets $A$ the difference set contains more elements than the sumset.   It had been a conjectured (cf. Croft~\cite{crof67}, Marica~\cite{mari69}, Nathanson~\cite{nath07a,nath07f}) that $|d(A)| \geq |s(A)|$ for every set $A$, but the set $A = \{0,2,3,4,7,11,12,14\}$ is a counterexample, since
\[
|A-A| = |d(A)|= 25
\]
and
\[
|A+A| = |s(A)|= 26.
\]
Sets with more sums than differences have also been studied by Hegarty~\cite{hega06} and Martin and O'Bryant~\cite{mart-obry06}.
This suggests the following problem.

\bp  \label{BLF:prob1}
Let $f(x_1,\ldots,x_n)$ and $g(x_1,\ldots,x_n)$ be polynomials with integer coefficients.  Determine if there exist finite sets $A, B, C$ of positive integers with $|C|>1$ such that 
\beq  \label{BLF:ABCineq}
\left\{\begin{array}{ll}
|f(A)|  >  |g(A)| \\
|f(B)|  <  |g(B)| \\
|f(C)|  =  |g(C)| .
\end{array}
\right.
\eeq
\ep

There is a stronger form of Problem~\ref{BLF:prob1}.
\bp  \label{BLF:prob4}
Let $f(x_1,\ldots,x_n)$ and $g(x_1,\ldots,x_n)$ be polynomials with integer coefficients.  Does there exist a sequence $\{A_i\}_{i=1}^{\infty}$ of finite sets of integers such that
\[
\lim_{i\rightarrow\infty} \frac{|f(A_i)|}{|g(A_i)|} = \infty ?
\]
Does there exist a sequence $\{C_i\}_{i=1}^{\infty}$ of finite sets of integers such that $\lim_{i\rightarrow\infty}|C_i| = \infty$ and 
$|f(C_i)| = |g(C_i)|$ for all $i$?
\ep

Linear polynomials constitute an important special case.

\bp  \label{BLF:prob2}
Let $f(x_1,\ldots,x_n) = u_1x_1 + \cdots + u_nx_n$ and $g(x_1,\ldots,x_n) = v_1x_1 + \cdots + v_nx_n$ be linear forms with integer coefficients.  Do there exist finite sets $A, B, C$ of integers with $|C|>1$ that satisfy~(\ref{BLF:ABCineq})?
\ep

The \emph{interval of integers} $[a,b]$ is the set of integers $\{a,a+1,a+2,\ldots,b\}.$  For any integer $u$ and sets $A$ and $B$ of integers, 
we define the \emph{dilation}
\[
u\ast A = \{ua:a\in A\}
\]
and the \emph{sumset} 
\[
A+B = \{a+b : a\in A \text{ and } b\in B\}.
\]
If $f(x_1,\ldots,x_n) = u_1x_1 + \cdots + u_nx_n$ is a linear form, then
\[
f(A) = u_1\ast A +\cdots + u_n\ast A.
\]
Sets of integers $A$ and $B$ are \emph{affinely equivalent} if there are rational numbers $u \neq 0$ and $v$ such that $B=u\ast A + \{v\}.$
In this case,
\[
f(B) = f(u\ast A + \{v\}) = u\ast f(A) + \{f(v,\ldots,v)\},
\]
hence
\[
|f(A)| = |f(B)|.
\]
Note that every two-element set is affinely equivalent to the set $\{0,1\},$
and that every set $A$ of integers with $1 < |A| < \infty$ is affinely equivalent to a set $A'$ such that $0 \in A'$ and $A'\setminus \{0\}$ is a set of relatively prime positive integers.
The following theorem implies that Problem~\ref{BLF:prob2} is equivalent to Problem~\ref{BLF:prob4} in the case of linear forms.

\bt
Let $f(x_1,\ldots,x_n)$ and $g(x_1,\ldots,x_n)$ be linear forms with integer coefficients, and let $A$ and $C$ be finite sets of  integers such that 
$|f(A)| > |g(A)|$, $|f(C)| = |g(C)|$, and $|C|>1$.  There exist  sequences $\{A_i\}_{i=1}^{\infty}$ 
and $\{C_i\}_{i=1}^{\infty}$ of finite sets of integers with 
$\lim_{i\rightarrow\infty}|C_i| = \infty $ such that
\[
\lim_{i\rightarrow\infty} \frac{|f(A_i)|}{|g(A_i)|} = \infty,
\]
and $|f(C_i)| = |g(C_i)|$ for all $i$.
\et

\begin{proof}
Let $f(x_1,\ldots,x_n)=\sum_{i=1}^n u_ix_i$ and $g(x_1,\ldots,x_n)=\sum_{i=1}^n v_ix_i$ be linear forms with integer coefficients, and let $A$ be a finite set of  integers.
We define 
\[
m_{f,g}(A) = \max\left( |s|: s \in A \cup f(A) \cup g(A)\right),
\]
and choose an integer
\[
M > 2m_{f,g}(A).
\]
Let
\[
A_M = A + M\ast A = \{a+Ma' : a,a' \in A\}.
\]
If $a_1,a'_1,a_2,a'_2 \in A$ and $a_1+Ma'_1 = a_2+ Ma'_2,$
then $a_1-a_2=M(a'_2-a'_1)$.  Since 
\[
M|a'_2-a'_1| = |a_1-a_2| \leq |a_1| + |a_2| \leq 2m_{f,g}(A) <M,
\]
it follows that $a'_1=a'_2$ and so  $a_1=a_2$ and
\[
|A_M| = |A|^2.
\]
The identity
\[
\sum_{i=1}^n u_i \left( a_i+Ma'_i \right)
= \sum_{i=1}^n u_i a_i + M \sum_{i=1}^n u_i a'_i
\]
implies that 
\[
f(A_M) = f(A)+M\ast f(A).
\]
If $s_1,s'_1,s_2, s'_2 \in f(A)$ and
\[
s_1+Ms'_1 = s_2+ Ms'_2,
\]
it again follows that $s_1=s_2,$ $s'_1=s'_2$, and 
\[
|f(A_M)| = |f(A)|^2.
\]
Similarly,
\[
|g(A_M)| = |g(A)|^2
\]
and so
\[
\frac{ |f(A_M)|}{|g(A_M)|} =  
\left(\frac{ |f(A)|}{|g(A)|} \right)^2.
\]
The Theorem follows by iterating this construction.
\end{proof}

In the case of binary linear forms, we write $f(x,y) = ux+vy$ instead 
of $f(x_1,x_2)=u_1x_1+u_2x_2.$     We are interested only in the cardinality of the image of $f(x,y)$ on a finite set $A$ of integers. We shall always assume that $uv \neq 0.$  If $(u,v)=d>1$ and $g(x,y) = (u/d)x+(v/d)y,$ then $|f(A)| = |g(A)|.$ Thus, we can assume that $(u,v)=1.$  Similarly, if $h(x,y) = vx+uy,$ then $|f(A)| = |h(A)|,$ and so we can assume that $|u| \geq |v|.$  Finally, if $\ell(x,y)  = -ux-vy,$ then $|f(A)| = |{\ell}(A)|,$ and we can assume that $u > 0.$  Therefore, it suffices to consider only binary linear forms $f(x,y)=ux+vy$ that have been normalized so that 
\[
u \geq |v| \geq 1 \text{ and }  (u,v)=1.
\]

\bp  \label{BLF:problem:binlinear}
Let $f(x,y) = u_1x + v_1y$ and $g(x,y) = u_2x + v_2y$ be normalized binary linear forms with nonzero integer coefficients $(u_1,v_1) \neq (u_2,v_2).$   Do there exist finite sets of integers $A$ and $B$ such that $|f(A)| > |g(A)|$ and  $|f(B)| < |g(B)|$?
\ep

In this paper we shall prove that the answer to the question in Problem~\ref{BLF:problem:binlinear} is ``yes''.

\section{Pairs of binary linear forms with $u_1,u_2 \geq 2$}  \label{BLF:section:u>2}

In this section we prove that if $f(x,y) = u_1x + v_1y$ and $g(x,y) = u_2x + v_2y$ are normalized binary linear forms with $u_1 \geq 2$,  $u_2 \geq 2$, and $(u_1,v_1) \neq (u_2,v_2),$ then there exist finite sets $A, B, C$ of integers such that   
$|f(A)| <  |g(A)|$, $|f(B)|> |g(B)|$, and $|f(C)| =  |g(C)|$.

\bt  \label{BLF:theorem:3set}
For $u > |v| \geq 1 \text{ and } (u,v)=1$, consider the normalized binary linear form
\[
f(x,y) = ux+vy.
\]
\benum
\item[(i)]
If $|A|=2,$ then $|f(A)|=4.$

\item[(ii)]
If  $u \geq 3$ and $|A| = 3$, then $|f(A)| = 8$ or 9, 
and $|f(A)| = 8$ if and only if $A$ is affinely equivalent to one of the two sets 
\[
\{0,|v|,u\} \text{ and } \{0,|v|,u+|v|\}.
\]

\item[(iii)]
If $u=2$ and $|A| = 3$, then $|f(A)| < 9$  if and only if $A$ is affinely equivalent to one of the two sets 
\[
\{0,1,2\}  \text{ and } \{0,1,3\}.
\]
Moreover, $|f(\{0,1,2\})| = 7$ and $|f(\{0,1,3\})|  = 8$.

\item[(iv)]
If $f(x,y) = ux+vy$, and $g(x,y) = ux-vy$, 
then $|f(A)| = |g(A)|$ for every set $A$ with $|A|=3.$
\eenum
\et

\begin{proof}
If $f(x,y) = ux+vy$ is a normalized binary linear form, then $f(\{0,1\}) = \{0,v,u,u+v\}$ and so $|f(\{0,1\})| = 4.$  Since every set $A$ with $|A|=2$ is affinely equivalent to $\{0,1\},$ it follows that if $|A|=2,$ then $|f(A)|=4.$
This proves~(i).

Let $|A| = 3.$  The set $A$ is affinely equivalent to a set $A'$ such that $\min(A')=0$ and $\gcd(A') = 1.$   If $|f(A)| \leq 8,$ then there exist $x_1,y_1, x_2,y_2 \in A'$ such that 
\[
ux_1+vy_1 = ux_2+vy_2 \text{ and } (x_1,y_1) \neq (x_2,y_2).
\]
It follows from~(i) that $|\{x_1,y_1,x_2,y_2\}| > 2$ and so
\[
\{x_1,y_1,x_2,y_2\} = \{a_1,a_2,a_3\} = A'.
\]
There are three possibilities: Either
\[
ua_1+va_2 = ua_1+va_3.
\]
or
\[
ua_1+va_2 = ua_3+va_1
\]
or
\[
ua_1+va_1 = ua_2+va_3.
\]

In the first case, $a_2=a_3,$ which is absurd.

In the second case, we have 
\[
u(a_1-a_3)=v(a_1-a_2).
\]
Since $(u,v)=1$, there exists an integer $r$ such that
\begin{align*}
a_1-a_2 & = ru  \\
a_1-a_3 & = rv \\
a_3-a_2 & = r(u-v).
\end{align*}
Since $0 \in A',$ it follows that $r$ divides each integer in $A'$, and so $r = \pm 1.$
If $a_1 = 0,$ then $r=-1$, $a_2=u$, $a_3=v=|v|,$ and $A' = \{0,|v|, u  \}.$ 
If $a_2 = 0,$ then $r=1$, $a_1=u$, and $a_3=u-v$.  If $v>0$, then $A' = \{0,u-|v|,u\}.$  If $v<0$, then $A' = \{0,u,u+|v|\}$.
If $a_3 = 0,$ then $a_2=-r(u-v)$ and so $r=-1$, $a_1=-v=|v|$,  $a_2=u-v=u+|v|$,and $A' = \{0, |v|, u+|v|  \}.$

In the third case, 
\[
u(a_1-a_2)=v(a_3-a_1)
\]
and there is an integer $r=\pm 1$ such that
\begin{align*}
a_1-a_2 & = rv  \\
a_3-a_1 & = ru \\
a_3-a_2 & = r(u+v).
\end{align*}
If $a_1 = 0,$ then $r=1$, $a_3 = u$, $a_2 = -v=|v|$, and $A' = \{0,|v|, u  \}.$
If $a_2 = 0,$ then $r=1$, $a_1=v=|v|$, $a_3=u+v=u+|v|$, and $A' = \{0,|v|,u+|v|  \}$.
If $a_3 = 0,$ then $r=-1$, $a_1=u$, and $a_2=u+v$.
If $v>0$, then $A' = \{0,u,u+|v| \}$.  If $v < 0$, then $A' = \{0,u-|v|,u \}$.

Since the two sets $\{0,|v|,u\}$ and $\{0,u-|v|,u\}$ are affinely equivalent, and the two sets $\{0,u,u+|v|\}$ and $\{0,|v|,u+|v|\}$ are also affinely equivalent, it follows that the sets $\{0,|v|,u\}$ and $ \{0,u,u+|v|\}$
are, up to affine equivalence, the only possible solutions of $|f(A)| \leq 8$ with $|A|=3.$   

We shall prove that if $u \geq 3$, then $|f(A)| = 8$ for both these sets.  
Let $v>0$ and $f(x,y) = ux+vy.$  If $A = \{0,v,u\}$, then 
\[
f(A) = \{ 0,v^2, uv, u^2, uv+v^2, 2uv, u^2 + v^2, uv + u^2\}.
\]
Since $2uv < u^2 + v^2$ for $v < u$, we have
\[
0 < v^2 < uv < uv+v^2 < 2uv < u^2+v^2 < uv + u^2
\]
and 
\[
uv < u^2 < u^2+v^2. 
\]
If $u^2 = 2uv,$ then $u=2v = 2$ since $(u,v)=1$.   If $u^2 = uv+v^2,$ then $u/v = (1+\sqrt{5})/2$, which is impossible since $u/v$ is rational.  Therefore, $|f(A)| = 8$.  

If $v>0$ and $A = \{0,v,u+v\}$, then  
\[
f(A) = \{ 0,v^2, uv, uv+v^2, 2uv+v^2, u^2+uv, u^2+uv+v^2,u^2+2uv+v^2\}.
\]
We have
\[
0 < v^2 < uv < uv+v^2 < 2uv + v^2 < u^2+uv+v^2 < u^2+2uv + u^2
\]
and 
\[
uv + v^2 < u^2+v^2 < u^2+uv+v^2. 
\]
If $u^2 +v^2= 2uv+v^2,$ then $u=2v = 2$,  which is false, and so $f(A)|=8.$
The case $u \geq 3$ and $v<0$ is similar.
This proves~(ii).

If $u = 2$, then $v = \pm 1$ and $f(x,y) = 2x \pm y.$  Up to affine equivalence, the sets $A$ with $|f(A)| \leq 8$ are $\{0,1,2\}$ and $\{0,1,3\}$.  For these sets, $|f(\{0,1,2\})| = 7$ and $|f(\{0,1,3\})| = 8$.
This proves~(iii).

To obtain~(iv), we observe that the binary linear forms $f(x,y) = ux+vy$ and $g(x,y) = ux-vy$ generate the same exceptional sets, and so $|f(A)| = |g(A)|$ whenever $|A|=3.$  This completes the proof.
\end{proof}

\bl   \label{BLF:lemma:homogen}
Let
\[
F(x,y) = c_0x^k + c_1x^{k-1}y + \cdots + c_{k-1}xy^{k-1} + c_ky^k
\]
be a nonzero homogeneous polynomial with integer coefficients.  Let $j$ be the largest integer such that $c_j \neq 0.$  If $u$ and $v$ are relatively prime nonzero integers such that $F(u,v)=0$, then $|c_j| \geq |u|.$  
\el

We call $c_j$ the \emph{last coefficient} in the polynomial $F(x,y)$.

\begin{proof}
If 
\[
F(u,v) = \sum_{i=0}^k c_iu^{k-i}v^i = u^{k-j} \sum_{i=0}^j c_iu^{j-i}v^i = 0
\]
then
\[
u\sum_{i=0}^{j-1} c_iu^{j-i-1}v^i = - c_jv^j
\]
and so $u$ divides $c_jv^j.$  Since $(u,v)=1$ and $c_j \neq 0$, it follows that $u$ divides $c_j$ and so  $|c_j| \geq |u|.$  
\end{proof}

\bt
Let 
\[
f(x,y) = u_1x+v_1y
\]
and 
\[
g(x,y) = u_2x+v_2y
\]
be normalized binary linear forms with 
\[
u_1\geq 2, u_2 \geq 2,
\]
and
\[
(u_1,|v_1|) \neq (u_2,|v_2|).
\]
There exist sets $A$ and $B$ with $|A| = |B| = 3$ such that
\[
|f(A)| < |g(A)| \text{ and } |f(B)| > |g(B)|.
\]
\et

\begin{proof}
If $u_1 < u_2$ and $u_2 \neq u_1+|v_1|,$ then the sets $A=\{0, |v_1|, u_1\}$ and $B = \{0,|v_2|,u_2 \}$ satisfy $|f(A)| = |g(B)| \leq 8$ 
and $|f(B)| = |g(A)| = 9$.
If $u_1 < u_2$ and $u_2 = u_1+|v_1|,$ then $u_2+|v_2| > u_1+|v_1|$ and the sets $A=\{0, |v_1|, u_1\}$ and $B = \{0,|v_2|,u_2+|v_2|\}$ satisfy $|f(A)| = |g(B)| = 8$ and $|f(B)| = |g(A)| = 9$.

If $u_1 = u_2$ and $|v_1| < |v_2|$, then sets $A=\{0,|v_1|,u_1+ |v_1|\}$ and $B = \{0,|v_2|,u_2+|v_2|\}$ satisfy $|f(A)| = |g(B)| = 8$ and $|f(B)| = |g(A)| = 9$.  This completes the proof.
\end{proof}

\bt
Let
\[
f(x,y) = ux+vy
\]
and 
\[
g(x,y) = ux-vy
\]
be normalized binary linear forms with $u > v \geq 1.$
For $u =2$, if
\[
A = \{ 0, 3,4,6\}
\]
and
\[
B = \{0,4,6,7\}
\]
then 
\[
|f(A)| = 13> 12 = |g(A)|
\] 
and
\[
|f(B)| = 13 < 14 = |g(B)|.
\] 
For $u \geq 3$, if
\[
A = \{ 0, u^2-v^2,u^2,u^2+uv\}
\]
and
\[
B = \{0,u^2-uv,u^2-v^2,u^2\}
\]
then 
\[
|f(A)| = 14 > 13 = |g(A)|
\] 
and
\[
|f(B)| = 13 < 14 = |g(B)|.
\] 
\et

\begin{proof}
For $u=2$ and sets $A = \{ 0, 3,4,6\}$ and $B = \{0,4,6,7\}$, we have
\begin{align*}
f(A) & = \{0,3,4,6,8,9,10,11,12,14,15,16,18\} \\
g(A) & = \{-6,-4,-3,0,2,3,4,5,6,8,9,12\} \\
f(B) & = \{0,4,6,7,8,12,14,15,16,18,19,20,21\} \\
g(B) & = \{-7,-6,-4,0,1,2,4,5,6,7,8,10,12,14\}
\end{align*}
with $|f(A)|=|f(B)|=13$, $|g(A)|=12$, and $|g(B)|=14.$

Let $u \geq 3$ and $A = \{ 0, u^2-v^2,u^2,u^2+uv\}$.  We list the elements of the set $f(A) = \{ux+vy : x,y \in A\}$ in the following table:

\begin{center}
\begin{tabular}{|c|cccc|}
\hline
$f(A)$& $0$ & $u^2-v^2$ & $u^2$ & $u^2+uv$  \\
\hline
$0$ & $0$ & $u^2v-v^3$ & $u^2v$ & $u^2v+uv^2$ \\
$u^2-v^2$ & $u^3-uv^2$ & $u^3+u^2v-uv^2-v^3$ & $u^3+u^2v-uv^2$ & \framebox{$u^3+u^2v$} \\
$u^2$ &  $u^3$ & $u^3+u^2v-v^3$ & \framebox{$u^3+u^2v$} & $u^3+u^2v+uv^2$   \\
$u^2+uv$ & \framebox{$u^3+u^2v$} & $u^3+2u^2v-v^3$ & $u^3+2u^2v$ & $u^3+2u^2v+uv^2$  \\
\hline
\end{tabular}
\end{center}

The number $u^3+u^2v$ occurs three times in this table, 
and so $|f(A)| \leq 14.$  By Lemma~\ref{BLF:lemma:homogen}, if two numbers in the table are equal for positive integers $u$ and $v$ with $u \geq 3$ and $(u,v)=1$,  then the difference of the two numbers is an expression of the form $F(u,v)$, where $F(x,y)$ is a homogeneous polynomial of degree 3 with last coefficient at least 3.  Since $1 \leq v < u$, the numbers in the table are increasing from left to right in each row and from top to bottom in each column.  The following 10  numbers in the set $f(A)$ are strictly increasing:
\begin{align*}
0 & < u^2v-v^3 < u^3-uv^2 < u^3+u^2v-uv^2-v^3  \\
& < u^3+u^2v-uv^2  < u^3+u^2v-v^3  < u^3+u^2v  \\
& < u^3+u^2v+uv^2 < u^3+2u^2v < u^3+2u^2+uv^2
\end{align*}
The four other numbers in $f(A)$ satisfy  
\[
u^2v < u^3 < u^3+2u^2v-v^3
\]
and 
\[
u^2v < u^2v+uv^2 < u^3+2u^2v-v^3.
\]
Comparing numbers among the three chains of inequalities, we see that there is no difference with last coefficient greater than 2, and so $|f(A)| = 14$.

Consider now the set $g(A) = \{ux-vy: x,y \in A \},$ whose elements we list in the following table:

\begin{center}
\begin{tabular}{|c|cccc|}
\hline
$g(A)$& $0$ & $u^2-v^2$ & $u^2$ & $u^2+uv$  \\
\hline
$0$ & $0$ & $-u^2v+v^3$ & $ -u^2v$ & $ -u^2v-uv^2$ \\
$u^2-v^2$ & \framebox{$u^3-uv^2$} & $u^3-u^2v-uv^2+v^3$ & \framebox{$u^3-u^2v-uv^2$} & $u^3-u^2v-2uv^2$ \\
$u^2$ &  \framebox{$u^3$} & $u^3-u^2v+v^3$ & $u^3-u^2v$ & \framebox{$u^3-u^2v-uv^2$}   \\
$u^2+uv$ & $u^3+u^2v$ & $u^3+v^3$ & \framebox{$u^3$} & \framebox{$u^3-uv^2$}  \\
\hline
\end{tabular}
\end{center}

The numbers $u^3, u^3-uv^2$, and $u^3-u^2v-uv^2$ occur two times in the table, and so $|g(A)| \leq 13.$  The numbers in the table are decreasing from left to right in each row and increasing from top to  bottom in each column.  The following 9 numbers in $g(A)$ form a strictly increasing sequence:
\begin{align*}
-u^2v-uv^2 & < -u^2v < -u^2v+v^3  < 0 < u^3-u^2v \\
& < u^3-u^2v+v^3  < u^3 < u^3 + v^3 < u^3+3u^2v.
\end{align*}
The other four numbers satisfy
\[
u^3-u^2v-2uv^2 < u^3-u^2v-uv^2 < u^3-u^2v-uv^2+v^3 < u^3-uv^2.
\]
Indeed, there is no pair of  expressions in the table whose difference has last coefficient greater than 2, and so $|g(A)|=13.$

Finally, we consider the sets $f(B)$ and $g(B)$:
\begin{center}
\begin{tabular}{|c|cccc|} 
\hline
$f(B)$& $0$ & $u^2-uv$ & $u^2-v^2$ & $u^2$  \\
\hline
$0$ & $0$ & $u^2v-uv^2$ & $u^2v-v^3$ & $u^2v$ \\
$u^2-uv$ & $u^3-u^2v$ & \framebox{$u^3-uv^2 $} & $u^3-v^3 $ & \framebox{$u^3$} \\
$u^2-v^2$ &  \framebox{$u^3-uv^2$} & $u^3+u^2v-2uv^2$ & $u^3+u^2v-uv^2-v^3$ & \framebox{$u^3+u^2v-uv^2$}   \\
$u^2$ & \framebox{$u^3$} & \framebox{$u^3+u^2v-uv^2$} & $u^3+u^2v-v^3$ & $u^3+u^2v$  \\
\hline
\end{tabular}
\end{center}

\begin{center}
\begin{tabular}{|c|cccc|} 
\hline
$g(B)$& $0$ & $u^2-uv$ & $u^2-v^2$ & $u^2$  \\
\hline
$0$ & $0$ & $-u^2v+uv^2$ & $-u^2v+v^3$ & $-u^2v$ \\
$u^2-uv$ & \framebox{$u^3-u^2v$} & $u^3-2u^2v+uv^2 $ & $u^3-2u^2v+v^3 $ & $u^3-2u^2v $ \\
$u^2-v^2$ &  $u^3-uv^2$ & \framebox{$u^3-u^2v$} & $u^3-u^2v-uv^2+v^3$ & $u^3-u^2v-uv^2$   \\
$u^2$ & $u^3$ & $u^3-u^2v+uv^2$ & $u^3-u^2v+v^3$ & \framebox{$u^3-u^2v$}  \\
\hline
\end{tabular}
\end{center}

In the table for $f(B)$, the numbers $u^3, u^3-uv^2,$ and $u^3+u^2v-uv^2$ occur twice, and so $|f(B)| \leq 13.$  
In the table for $g(B)$, the number $u^3-u^2v$ occurs three times, so $|g(B)| \leq 14.$  In neither table is there a pair of numbers whose difference has last coefficient greater than 2, and so $|f(B)| = 13$  
and $|g(B)| = 14.$ 
This completes the proof.
\end{proof}

\bt
Let $u$ and $v$ be relatively prime positive integers with $u > v,$ and consider the forms
\[
f = ux+vy
\]
and
\[
g = ux-vy.
\]
If $A$ is an arithmetic progression of length $t \leq u,$ then  $|f(A)| = |g(A)| = t^2.$
\et

\begin{proof}
Since an arithmetic progression of length $t$ is affinely equivalent to the interval $[0,t-1]$, it suffices to consider the sets $A_t = [0,t-1]$ for $t = 1,\ldots,u.$

If $x_1,x_2,y_1,y_2 \in A_t$ and $ux_1+vy_1 = ux_2+vy_2,$ then 
$u(x_1-x_2) = v(y_2-y_1).$  Since $(u,v)=1,$ it follows that $u$ divides $y_2-y_1.$  Since $|y_2-y_1| < t\leq u,$ it follows that $y_1=y_2,$  
and so $x_1=x_2$.  Thus, every element in $f(A_t)$ has a unique representation in the form $ux+vy,$ and $|f(A_t)|=t^2.$  The proof that $|g(A_t)|=t^2$ is similar.
\end{proof}

\section{The pair of linear forms $ux+vy$ and $x-y$}  \label{BLF:section:f-d}
Let $u$ and $v$ be relatively prime positive integers with $u > v$, and consider the linear forms 
\[
f(x,y) = ux+vy \text{ and } d(x,y) = x-y.
\]
\bt
Let 
\[
A = \{ 0, v^3, v^3+v^2u,v^3+v^2u+vu^2, v^3+v^2u+vu^2+u^3\}.
\]
Then 
\[
f(A) \leq 19 \text{ and } d(A) = 21.
\]
\et

\begin{proof}
Let
\begin{align*}
a_0 & = 0 \\
a_1 & = v^3 \\
a_2 & = v^3+v^2u \\
a_3 & = v^3 + v^2u + vu^2 \\
a_4 & = v^3 + v^2u + vu^2 + u^3.
\end{align*}
Then $A = \{a_0, a_1,a_2,a_3,a_4\},$ and
\[
a_0 < a_1 < a_2 < a_3 < a_4.
\]
Since $|A|=5,$ we have $|f(A)| \leq 25$ and $|d(A)| \leq 21.$
To show that $|f(A)| \leq 19$, it suffices to give six different integers, each of which has two distinct representations  in $f(A)$.  Here they are:
\begin{align*}
n_1 & = ua_1+va_1 = ua_0 + va_2 \\
n_2 & = ua_2+va_1 = ua_0 + va_3 \\
n_3 & = ua_2+va_2 = ua_1 + va_3 \\
n_4 & = ua_3+va_1 = ua_0 + va_4 \\
n_5 & = ua_3+va_2 = ua_1 + va_4 \\
n_6 & = ua_3+va_3 = ua_2 + va_4. 
\end{align*}
A straightforward calculation shows that
\[
n_1 < n_2 < n_3 < n_4 < n_5 < n_6
\]
and so $|f(A)| \leq 19.$

Next we prove that $|d(A)| = 21$.  Let $D = \{a_j-a_i: 0 \leq i < j \leq 4\}$.
It suffices to prove that $|D| = 10$.  If $v=1$, then $A = \{1,1+u,1+u+u^2,1+u+u^2+u^3\}$ is a Sidon set and $|D|=10.$  

Suppose that $v \geq 2$. 
Since $u>v$ and $a_0 = 0$, we have
\[
a_1 < a_2-a_1 < a_3-a_2 <  a_4 - a_3 < a_4-a_2 < a_4 - a_1 < a_4
\]
and
\[
a_2 < a_3 - a_1 < a_3.
\]
Let
\[
D_1 = \{ a_1 , a_2-a_1 , a_3-a_2 ,  a_4 - a_3 , a_4-a_2 , a_4 - a_1 , a_4 \} 
\]
and 
\[
D_2 = \{a_2 , a_3 - a_1 , a_3\} .
\]
We must show that $D_1\cap D_2 = \emptyset$.
There are three cases.

Case ($a_2$):  Since
\[
a_2 -a_1 < a_2 < a_4-a_2,
\]
it follows that if $a_2 \in D_1$, then $a_2 = a_3-a_2$ or $a_2 = a_4-a_3$.
If $a_2 = a_3-a_2 $, then $v^2(v+u) = v^3+v^2u = vu^2$, and  $v(v+u) = u^2$.  Since $(u,v)=1$, it follows that $v=1$.
If $a_2 = a_4-a_3$, then $v^2(v+u) = v^3+v^2u =u^3$ and $v=1$.

Case ($a_3-a_1$):  Since
\[
a_3-a_2 < a_3 - a_1 < a_4-a_2,
\]
it follows that if $a_3-a_1 \in D_1$, then $a_3-a_1 = a_4-a_3$.
If $a_3 - a_1 = a_4-a_3$, then $vu(v+u) = v^2 u +v u^2 = u^3$, and so $v(v+u) = u^2$.  This implies that $v=1$.

Case ($a_3$):   Since
\[
a_3 -a_2 < a_3 < a_4-a_1,
\]
it follows that $a_3 = a_4-a_3$ or $a_3 = a_4-a_2$.
If $a_3 = a_4-a_3 $, then $v(v^2+vu+u^2) = v^3+v^2u+vu^2 = u^3$, and $v=1$.  If $a_3 = a_4-a_2 $, then $v^3+v^2u+vu^2 = vu^2 + u^3$, 
hence $v^2(v+u) = u^3$ and $v=1$.  This completes the proof.
\end{proof}

\section{The pair of linear forms $ux+vy$ and $x+y$} 

In Section~\ref{BLF:section:u>2} we solved Problem~\ref{BLF:problem:binlinear} for pairs of normalized binary linear forms $f(x,y) = u_1x + v_1y$ and $g(x,y) = u_2x + v_2y$ with $u_1,u_2 \geq 2.$    In Section~\ref{BLF:section:f-d} we solved the case $f(x,y) = ux + vy$ with $u \geq 2$ and $d(x,y) = x-y$.  The remaining case is $f(x,y) = ux+vy$ with $u \geq 2$ and $s(x,y) = x +y$.   

For example, consider the form $f(x,y) = 2x+y$.  We have
\[
4 = |f(\{0,1\})|  >   |s(\{0,1\})| = 3.
\]
We shall construct a set $A$ with $ |f(A)|  <   |s(A)|.$
Start by defining the four sets:
\begin{align*}
	R_{13} & = \{0,1,6,7,9,11\}  \\
	R_{15} & = \{0,1,5,6,10,11,13 \} \\
	R_{16} & = \{0,1,3,5,7,9,11,13,15 \} \\
	R_{19} & = \{0,1,11,12,14,16,18\}.
\end{align*}
Note that 
\[
13\cdot 15 \cdot 16 \cdot 19= 59280
\]
and
\[
|R_{13}|\cdot |R_{15}| \cdot |R_{16}| \cdot |R_{19}|= 6\cdot 7\cdot 9\cdot 7 = 2646.
\]
Let $x \mod m$ denote the least nonnegative integer that is congruent to $x$ modulo $m$.  We define
\[
A =	\{x \in [1, 59280]: x \mod m \in R_m  \text{ for all } m \in \{13,15,16,19\}\}.
\]
The set $A$ contains  2646 elements.
By direct calculation, we have $|f(A)| = 108014$ and $|s(A)|= 114575.$

The linear form $f(x,y) = 2x+y$ is a special case.  In general, we do not have an algorithm to construct finite sets $A$ of integers such that $|f(A)| < |s(A)|$  for an arbitrary normalized bilinear form $f(x,y) = ux+vy$ with $u \geq 2.$   However, such sets do exist.
In the following sections we shall show that, associated to the form $f(x,y) = ux+vy$, there is an infinite set $M$ of positive integers with the property that, for each $m \in M,$ there is a set of congruence classes $R_m \subseteq \Z/m\Z$ such that $s(R_m) = \Z/m\Z$ and $f(R_m) \subsetneqq \Z/m\Z.$   From the sets $R_m$ we construct a finite set $A$ of nonnegative integers such that  $|f(A)| < |s(A)|.$   Thus, we combine local solutions of the inequality to construct a global solution.

\section{A local to global criterion for pairs of linear forms in $n$ variables}

\bl      \label{BLF:lemma:crt}
Let $f(x_1,\ldots,x_n)$ be a polynomial with integer coefficients.
Let $m_1,\ldots, m_r$ be pairwise relatively prime positive integers, and $m = m_1\cdots m_r.$
Let $R_{m_i}$ be a set of congruence classes modulo $m_i$ for $i = 1,\ldots,r.$ 
Let $R_m$ be the set of all congruence classes $a+m\Z$ such that $a+m_i\Z \in R_{m_i}$ for $i = 1,\ldots,r$.  Then
\[
|R_m| = \prod_{i=1}^r |R_{m_i}|
\]
and
\[
|f(R_m)| = \prod_{i=1}^r |f(R_{m_i})|.
\]
\el

\begin{proof}
This follows from the Chinese remainder theorem.
\end{proof}

\bl     \label{BLF:lemma:rectify}
Let $f(x_1,\ldots,x_n) = u_1x_1+ u_2x_2+\cdots u_n x_n$ be a linear form with integer coefficients, and let 
\[
h_f = |u_1|+|u_2|+\cdots + |u_n|.
\]
Let $R_m$ be  a set of congruence classes in $\Z/m\Z$, and let $A$ be the set of integers that consists of the least nonnegative element of each congruence class in $R_m.$  Then 
\beq  \label{BLF:ineq1}
|f(R_m)| \leq |f(A)| \leq 2h_f |f(R_m)|.
\eeq
\el

\begin{proof}
The triangle inequality implies that
\[
|f(a_1,\ldots,a_n)| \leq h_f \max(|a_i| : i=1,\ldots,n)
\]
for all integers $a_1,\ldots,a_n$.
Since $A \subseteq [0,m-1]$, it follows that
$f(A) \subseteq [-h_f (m-1),h_f (m-1)]$.	
If $a \in f(A),$ then $a+m\Z \in f(R_m).$  The lower bound in~(\ref{BLF:ineq1}) follows from the fact that $f(A)$ contains at least one element of every congruence class in $f(R_m)$.
The upper bound in~(\ref{BLF:ineq1}) follows from the fact that the interval $[-h_f (m-1),h_f(m-1)]$ contains at most $2h_f$ members of any congruence class modulo $m$.
\end{proof}

\bt    \label{BLF:theorem:rectify-main}
Let $f(x_1,\ldots,x_n) = u_1x_1+ \cdots u_n x_n$ and $g(x_1,\ldots,x_n) = v_1x_1+ \cdots v_n x_n$ be binary linear forms.
Let $M = \{m_i\}_{i=1}^{\infty}$ be a set of pairwise relatively prime integers  such that $m_i \geq 2$ for all $m_i \in M$.  If for every $m_i \in M$ there exists a nonempty set $R_{m_i}$ of congruence classes in $\Z/m_i\Z$ such that
\beq   \label{BLS:infprod}
\prod_{i=1}^{\infty} \frac{|f(R_{m_i})|}{|g(R_{m_i})|} = 0,
\eeq
then there is a finite set $A$ of integers such that 
\[
|f(A)| < |g(A)|.
\]
\et

\begin{proof}
Let 
\[
h_f = |u_1|+|u_2|+\cdots + |u_n|.
\]
Since the infinite product~(\ref{BLS:infprod}) diverges to 0, there is an integer $r$ such that 
\[
\prod_{i=1}^{r} \frac{|f(R_{m_i})|}{|g(R_{m_i})|} < \frac{1}{2h_f}.
\]
Let $m = m_1\cdots m_r$ and let  $R_m$ be the set of all congruence classes $a+m\Z$ such that $a+m_i\Z \in R_{m_i}$ for $i = 1,\ldots,r$.  
By Lemma~\ref{BLF:lemma:crt},
\[
 \frac{|f(R_m)|}{|g(R_m)|} = \prod_{i=1}^{r} \frac{|f(R_{m_i})|}{|g(R_{m_i})|}.
\]
Let $A$ be the set of integers that consists of the least nonnegative element in each congruence class in $R_m.$  By Lemma~\ref{BLF:lemma:rectify} we have
\[
|f(A)|  \leq 2h_f |f(R_m)|  < |g(R_m)|  \leq |g(A)|.
\]
This completes the proof.
\end{proof}

\bt    \label{BLF:theorem:rectifycong}
Let $f(x_1,\ldots,x_n) = u_1x_1+ \cdots u_n x_n$ and $g(x_1,\ldots,x_n) = v_1x_1+ \cdots v_n x_n$ be binary linear forms.
Let $M = \{m_i\}_{i=1}^{\infty}$ be a set of pairwise relatively prime positive integers  such that $m_i \geq 2$ for all $m_i \in M$ and 
\[
\sum_{i=1}^{\infty} \frac{1}{m_i} = \infty.
\]
If for every $m_i \in M$ there exists a nonempty set $R_{m_i}$ of congruence classes in $\Z/m_i\Z$ such that
\[
f(R_{m_i}) \neq \Z/m_i\Z
\]
and
\[
g(R_{m_i}) = \Z/m_i\Z,
\]
then there is a finite set $A$ of integers such that 
\[
|f(A)| < |g(A)|.
\]
\et

\begin{proof}
Since $|f(R_{m_i})| \leq m_i-1$ and $|g(R_{m_i})| = m_i$ for all $m_i \in M$, we have
\[
\frac{|f(R_{m_i})|}{|g(R_{m_i})|} \leq 1- \frac{1}{m_i}.
\]
The divergence of the infinite series $\sum_{i=1}^{\infty} m_i^{-1}$ implies that
\[
\prod_{i=1}^{\infty} \frac{|f(R_{m_i})|}{|g(R_{m_i})|} = \prod_{i=1}^{\infty}\left( 1- \frac{1}{m_i}\right) = 0,
\]
and the result follows immediately from Theorem~\ref{BLF:theorem:rectify-main}.  
\end{proof}

We can restate Theorem~\ref{BLF:theorem:rectifycong} as follows.

\bt   \label{BLF:theorem:rectifyint}
Let $f(x_1,\ldots,x_n) = u_1x_1+ \cdots u_n x_n$ and $g(x_1,\ldots,x_n) = v_1x_1+ \cdots v_n x_n$ be binary linear forms.
Let $M = \{m_i\}_{i=1}^{\infty}$ be a set of pairwise relatively prime positive integers  such that $m_i \geq 2$ for all $m_i \in M$ and 
\[
\sum_{i=1}^{\infty} \frac{1}{m_i} = \infty.
\]
If for every $m_i \in M$ there exists an integer $q_{m_i}$ and a finite set $A_{m_i}$ of integers such that
\benum
\item[(i)]  $f(a_1,\ldots,a_n) \not\equiv q_{m_i} \pmod{m_i}$ for all $a_1,\ldots,a_n \in A_{m_i},$ and 
\item[(ii)]
for every integer $q$ the congruence $g(a_1, \ldots, a_n) \equiv q \pmod{m}$ is solvable with $a_1,\ldots,a_n \in A_{m_i},$ 
\eenum
then there is a finite set $A$ of integers such that 
\[
|f(A)| < |g(A)|.
\]
\et

\section{An application of quadratic reciprocity}
\bt    \label{BLF:theorem:sumquad}
Let $p$ be a prime number such that $p \equiv 1 \pmod{4}$ and $p>5,$
and let
\[
R_p = \{k^2+p\Z : k=1,2,\ldots,(p-1)/2\}
\]
be the set of quadratic residues modulo $p$.
Let $s(x,y) = x+y$ and $d(x,y) = x-y.$  Then
\[
s(R_p) = d(R_p) = \Z/p\Z.
\]
\et

\begin{proof}
Let $a \in \Z.$  
In the finite field $\F_p = \Z/p\Z$, consider the sets 
\[
X = \{x^2+p\Z : x=0,1,\ldots,(p-1)/2\}
\]
and
\[
Y = \{a-y^2+p\Z:y=0,1,\ldots,(p-1)/2\}.
\]
Since $|X|=|Y|=(p+1)/2,$ it follows from the pigeonhole principle that $X\cap Y \neq \emptyset,$ and so there exist integers $x$ and $y$ such that $x^2 \equiv a -y^2 \pmod{p},$ that is, 
\beq   \label{BLF:cong1}
x^2 + y^2 \equiv a \pmod {p}.
\eeq
We must show that if $p \equiv 1 \pmod{4}$ and $p > 5,$ then this congruence can always be solved with integers not divisible by $p$.

Let $(a|p)$ denote the Legendre symbol modulo $p.$  For primes $p \equiv 1 \pmod{4},$ we have $(-1|p) = 1,$ and there is an integer $w$ such that $(w,p)=1$ and 
\[
1^2 + w^2 \equiv 0 \pmod{p}.
\]
Thus, we only have to consider congruences of the form~(\ref{BLF:cong1}) with $a \not\equiv 0 \pmod{p}.$  Moreover, we can assume that $a$ is a quadratic residue modulo $p$, since, if $x \equiv 0 \pmod{p},$ then $y^2 \equiv a\pmod{p}.$   

We begin with the case $a \equiv 1 \pmod{p}.$  At least one of the integers $2,3,6$ is a quadratic  residue modulo $p$, since if $(2|p)=(3|p)=-1,$ then $(6|p)=(2|p)(3|p)=1.$  There are three cases.
\benum
\item[(i)]
If $(2|p) = 1,$ then there is an integer $r$ such that 
$
r^2 \equiv 2 \pmod{p}$ and so
\[
3^2 + (2wr)^2 \equiv 9-8 \equiv 1 \pmod{p}.
\]
\item[(ii)]
If $(3|p) = 1,$ then there is an integer $s$ such that 
$
s^2 \equiv 3 \pmod{p}
$
and so
\[
2^2 + (ws)^2 \equiv 4-3 \equiv 1 \pmod{p}.
\]

\item[(iii)]
If $(6|p) = 1,$ then there is an integer $t$ such that 
$
t^2 \equiv 6 \pmod{p}
$
and so
\[
5^2 + (2wt)^2 \equiv 25-24 \equiv 1 \pmod{p}.
\]
\eenum
Since $p > 5,$, there exist integers $c$ and $d$ not divisible by $p$ such that
\[
c^2 + d^2 \equiv 1 \pmod{p}.
\]
If $a \equiv y^2\pmod{p},$ then
\[
(cy)^2 + (dy)^2 \equiv y^2 \equiv a \pmod{p}
\]
and
\[
(cy)^2 - (dwy)^2 \equiv y^2 \equiv a \pmod{p}.
\]
This completes the proof.
\end{proof}

\bt   \label{BLF:theorem:badp}
Consider the binary linear form
\[
f(x,y)=ux+vy
\]
where $u$ and $v$ are integers not divisible by $p$.
Let $p$ be a prime number and let $R_p$ be the set of quadratic residues modulo $p$.   Then $p\Z \in f(R_p)$ if and only if $-uv$ is a quadratic residue modulo $p$,
\et

\begin{proof}
If $p\Z \in f(R_p),$ then there are integers $k_1$ and $k_2$ not divisible by $p$ such that 
\[
uk_1^2 +vk_2^2\equiv 0\pmod{p}.
\]
Then
\[
uvk_1^2 + (vk_2)^2  \equiv 0\pmod{p}
\]
and so
\[
-uv \equiv  \left(vk_2k_1^{-1}\right)^2 \pmod{p}
\]
that is, $-uv$ is a quadratic residue modulo $p$.

Conversely, if $-uv$ is a quadratic residue modulo $p,$ then there is an integer $z \not\equiv 0 \pmod{p}$ such that $-uv\equiv z^2\pmod{p}$ and so 
\[
f(v^2,z^2) = uv^2+vz^2 \equiv 0 \pmod{p}.
\]
Thus, $p\Z = f(v^2+p\Z,z^2+p\Z) \in f(R_p).$
\end{proof}

\bt   \label{BLF:theorem:uvnonsquare}
Let
\[
f(x,y)=ux+vy
\]
be a normalized bilinear form such that $|uv|$ is not a perfect square.
Let
\[
s(x,y) = x+y
\]
and
\[
d(x,y) = x-y.
\]
There exist finite sets $A$ and $A'$ of integers such that
\[
|f(A)| < |s(A)| \text{ and } |f(A')| < |d(A')|.
\]
\et

\begin{proof}
Since $f(x,y)$ is normalized and $|uv|$ is not a square, we can write
\[
|uv| = w^2 2^{\varepsilon}  \prod_{j=1}^t q_j 
\]
where $w$ is a positive integer, $\varepsilon \in \{0,1\}$,  $q_1,\ldots,q_t$ are distinct odd primes, and $q_i \neq p$ for $i=1,\ldots,p.$
For $p \equiv 1 \pmod{4},$  we have
\begin{align*}
\left(\frac{-uv}{p}\right) & = \left(\frac{\pm w^2 2^{\varepsilon}\prod_{j=1}^t q_j}{p}\right)  = \left(\frac{2}{p} \right)^{\varepsilon} \prod_{j=1}^t \left(\frac{q_j}{p}\right) = \left(\frac{2}{p} \right)^{\varepsilon}\prod_{j=1}^t \left(\frac{p}{q_j}\right).
\end{align*}
If $\varepsilon=1,$ we choose $p$ so that $p \equiv 5\pmod{8}$ and $p\equiv 1\pmod{q_j}$ for $j=1,\ldots,t.$  Then  $(-uv|p)=-1.$
If $\varepsilon=0,$ then $t \geq 1$ and we choose $p$ so that $p \equiv 1\pmod{4}$, $p\equiv 1\pmod{q_j}$ for $j=2,\ldots,t$, and $(p|q_1)=-1.$  Again, $(-uv|p)=-1.$
In both cases, there is at least one infinite arithmetic progression $P(u,v)$ such that if $p$ is a prime 
and $p \in P(u,v),$ then $p\equiv 1 \pmod{4}$ and $(-uv|p)=-1.$  By Dirichlet's theorem, the arithmetic progression $P(u,v) $ contains infinitely many primes and 
\[
\sum_{\substack{p\in P(u,v)\\ p \text{ prime}} } \frac{1}{p} = \infty.
\]
By Theorem~\ref{BLF:theorem:sumquad}   and Theorem~\ref{BLF:theorem:badp}, for each prime $p \in P(u,v),$ the set
\[
R_p=\{k^2+p\Z: k=1,2,\ldots,(p-1)/2\}
\]
satisfies
\[
f(R_p) \neq \Z/p\Z
\]
and
\[
s(R_p) = d(R_p) = \Z/p\Z.
\]
The result now follows from Theorem~\ref{BLF:theorem:rectifycong}.
\end{proof}

\section{An exponential sum}
After Theorem~\ref{BLF:theorem:uvnonsquare}, we are left to consider only normalized bilinear forms $f(x,y) = ux+vy$ such that $|uv|$ is a square.  Since $u$ and $v$ are relatively prime, if follows that there are positive integers $U$ and $V$ not divisible by $p$ such that $u=U^2$ and $v=\pm V^2$. 
 Let $R_p$ denote the set of quadratic residues modulo the prime $p$.  If $f(x,y) = U^2x+V^2y,$ then $f(R_p) = s(R_p)$.   If $f(x,y) = U^2x-V^2y,$ then $f(R_p) = d(R_p)$.  This suggests that considering only squares mod $p$ will not suffice to resolve Problem~\ref{BLF:problem:binlinear} in this remaining case.  We shall generalize our method to $k$th powers.  We begin by applying elementary harmonic analysis on finite fields to binary linear forms.  A general reference is Nathanson~\cite[chapter 4]{nath00aa}.

Let $p$ be a prime number and $\F_p = \Z/p\Z$ the field of congruence classes modulo $p$.  We denote the multiplicative group of the field by $\F_p^{\times}$ and define 
\[
e_p(t) = e^{2\pi i t/p}.
\]
For all integers $a, b,$ and $t$ we have
\[
\sum_{t=0}^{p-1} e_p((a-b)t) = 
\begin{cases}
0 & \text{ if $a \not\equiv b \pmod{p}$} \\
p & \text{ if $a \equiv b \pmod{p}$.}
\end{cases}
\]
If $x$ is a congruence class modulo $p$, that is, if $x = t + p\Z$ for some integer $t$, then we define $e_p(x) = e_p(t)$.  This function is well-defined on $\F_p.$ 

Let $\gamma$ be a complex-valued function on \Fp.
We define the \emph{Fourier transform} $\hat{\gamma}: \F_p \rightarrow \C$ by 
\[
\hat{\gamma}(x) = \sum_{y\in \Fp} \gamma(y)e_p(-xy).
\] 
We have
\[
\hat{\gamma}(0) = \sum_{y\in \Fp} \gamma(y)
\]
and Plancherel's formula~\cite[Theorem 4.9]{nath00aa}
\beq  \label{BLF:plancherel}
\sum_{x\in \Fp}|\hat{\gamma}(x)|^2 = p\sum_{x\in \Fp} |\gamma(x)|^2 .  \eeq

Let $H$ be a subset of \Fp\ of cardinality $n$.  We also use $H$ to denote the characteristic function of $H$, that is, $H:\Fp \rightarrow \C$ is the function defined by
\[
H(x) = \begin{cases}
1 & \text{ if $x\in H$} \\
0 & \text{ if $x\notin H$.}
\end{cases}
\]
Then
\[
\hat{H}(x) = \sum_{y\in \Fp} H(y)e_p(-xy) = \sum_{h\in H}e_p(-xh)
\]
and
\[
\hat{H}(0)=\card(H)=n.
\]
Applying Plancherel's formula to the function $H$, we obtain
\[
\sum_{x\in \F_p} |\hat{H}(x)|^2 = p\sum_{x\in \F_p} |H(x)|^2= p\card(H) = pn
\]
and so
\beq  \label{BLF:mordell1}
\sum_{x\in \F_p^{\times}} |\hat{H}(x)|^2 = \sum_{x\in \F_p} |\hat{H}(x)|^2 - |\hat{H}(0)|^2 =(p-n)n.
\eeq

\bt  \label{BLF:theorem:mordell}
Let $f(x,y) = ux+vy$, where $u$ and $v$ are integers not divisible by $p$.  Let $H$ be a subgroup of order $n \geq 2$ of the multiplicative group $\F_p^{\times}$  and let 
\[
k = [\F_p^{\times}:H] = \frac{p-1}{n}.
\]
If $p > k^4,$ then $\Fp^{\times} \subseteq f(H)$, that is, every element of $\F_p^{\times}$ can be represented in the form $f(h_1,h_2)$ for some $h_1,h_2\in H$.
\et

\begin{proof}
Define the \emph{representation function} $r:\Fp \rightarrow \N_0$ as follows.
For every $x \in \F_p,$ let $r(x)$ denote the number of ordered pairs $(h_1,h_2) \in H \times H$ such that $f(h_1,h_2) = x$.   Then
\[
\sum_{x\in \F_p} r(x) = |H|^2 = n^2
\]
and
\begin{align*}
\hat{r}(x) 
& = \sum_{y\in \F_p} r(y)e_p(-xy) \\
& = \sum_{y\in \F_p} 
\left( \sum_{ \substack{h_1,h_2 \in H\\ uh_1+vh_2 = y} } 1\right) 
e_p(-xy) \\
& = \sum_{h_1\in H}\sum_{h_2\in H} e_p(-(uh_1+vh_2)x) \\
& = \sum_{h_1\in H} e_p(-uxh_1) \sum_{h_2\in H} e_p(-vxh_2) \\
& = \hat{H}(ux) \hat{H}(vx).
\end{align*}

Applying Plancherel's formula~\eqref{BLF:plancherel} to the function $r(x)$, we obtain
\begin{align*}
\sum_{x\in \F_p^{\times}} |\hat{H}(ux) & \hat{H}(vx)|^2
  =  \sum_{x\in \F_p} |\hat{H}(ux)\hat{H}(vx)|^2 - |\hat{H}(0|^4  \\
& = \sum_{x \in \F_p} \hat{r}(x)^2 - n^4 = p\sum_{x \in \F_p} r(x)^2 - n^4  \\
& = p\sum_{x \in \F_p} \left( r(x) - \frac{n^2}{p}\right)^2 .
\end{align*}

Let $x,x' \in \F_p^{\times}$.  If $x$ and $x'$ belong to the same coset of $\F_p/H$, then there exists $h'\in H$ such that $x = x'h'.$  It follows that
\[
\hat{H}(x) = \sum_{h\in H} e_p(-xh) 
 = \sum_{h\in H} e_p(-x'h'h) 
 = \sum_{h\in H} e_p(-x'h) 
=\hat{H}(x')
\]
and so the Fourier transform $\hat{H}(x)$ is constant on the cosets of $\F_p^{\times}/H$.
Similarly,  $uh_1+vh_2 = x'$  if and only if $u h_1h'+vh_2h' = x'h' = x$, and so $r(x)=r(x')$, that is, the representation function $r(x)$ is also constant on the cosets of $\F_p^{\times}/H$.

Let $\{x_1,\ldots,x_k\} \subseteq \F_p^{\times}$ be a set of coset representatives of $H$, that is, 
\[
\F_p^{\times}/H = \{x_1H,\ldots,x_kH\}.
\]
Applying~\eqref{BLF:mordell1}, we obtain
\[
(p-n)n = \sum_{x\in \F_p^{\times}} |\hat{H}(x)|^2 = \sum_{i=1}^k \sum_{x\in x_iH} |\hat{H}(x)|^2 = n\sum_{i=1}^k |\hat{H}(x_i)|^2 
\]
and so
\[
\sum_{i=1}^k |\hat{H}(x_i)|^2 = p-n.
\]
For every $x \in \F_p^{\times}$, there is an integer $j \in \{1,\ldots,k\}$ such that $x\in x_j H$ and $\hat{H}(x) = \hat{H}(x_j)$.  It follows that
\[
|\hat{H}(x)|^2 = |\hat{H}(x_j)|^2 \leq \sum_{i=1}^k |\hat{H}(x_i)|^2 = p-n.
\]
Since $(u,p) = (v,p)=1$, we have
\begin{align*}
p\sum_{x\in \F_p} \left( r(x) - \frac{n^2}{p}\right)^2
& = \sum_{x\in \F_p^{\times}} |\hat{H}(ux)\hat{H}(vx)|^2 \\
& \leq (p-n)\sum_{x\in \F_p^{\times}} |\hat{H}(ux)|^2 \\
& = (p-n)\sum_{x\in \F_p^{\times}} |\hat{H}(x)|^2  \\
& = (p-n)^2n.
\end{align*}
Since the representation function $r(x)$ is constant on cosets of $H$, we have
\begin{align*}
\sum_{i=1}^k \left( r(x_i) - \frac{n^2}{p}\right)^2
& = \frac{1}{n} \sum_{i=1}^k \sum_{x\in x_iH} \left( r(x) - \frac{n^2}{p}\right)^2 \\
& = \frac{1}{n}\sum_{x\in \F_p} \left( r(x) - \frac{n^2}{p}\right)^2 \\
& \leq \frac{(p-n)^2}{p}.
\end{align*}
For every $x \in \F_p^{\times}$, we have $x \in x_j H$ for some $j$,  and so 
\[
\left( r(x) - \frac{n^2}{p}\right)^2  =
\left( r(x_j) - \frac{n^2}{p}\right)^2 
\leq \frac{(p-n)^2}{p}.
\]
If $r(x)=0$ for some $x \in \F_p^{\times}$, then
\[
\frac{n^4}{p^2} \leq \frac{(p-n)^2}{p}.
\]
Since $|H|= n \geq 2$ and $p(p-n) \leq p(p-2) < (p-1)^2$, we have
\[
\left(\frac{p-1}{k} \right)^4 = n^4 \leq p(p-n)^2 < (p-1)^3
\]
and
\[
p \leq  k^4.
\]
This proves that if $p > k^4$, then $r(x) \geq 1$ for all $x \in \F_p^{\times}$ and so $\Fp^{\times} \subseteq f(H).$
\end{proof}

A finite cyclic group $G$ of order $N$ has a unique subgroup $H$ of order $n$ for every positive divisor $n$ of $N$.  If $k=[G:H] = N/n$ and if $g$ is a generator of $G$, then $H = \{g^{ik}: i=0,1,\ldots,n-1\} = \{x^k:x\in G\}$ is the set of $k$th powers in $G$.  

Let $H$ be a subgroup order $n$ of $\Fp^{\times}$ and let $k = [\Fp^{\times}:H] = (p-1)/n.$  Since the multiplicative group of a finite field is cyclic, it follows that $H$ is the subgroup of $k$th powers mod $p$, that is, $H = \{x^k:x\in \Fp^{\times}\}$.
We can restate Theorem~\ref{BLF:theorem:mordell} as follows.

\bt  \label{BLF:theorem:mordell-k}
Let $f(x,y) = ux+vy$ be a binary linear form with nonzero integral coefficients $u$ and $v$.  For $k\geq 2$, let $p$ be a prime number such that $p\equiv 1\pmod{k}$, $p > k^4$, and $(p,uv)=1$.  
If $H$ is the set of $k$th powers in $\Fp^{\times}$, then
\[
\Fp^{\times} \subseteq f(H).
\]
\et

We shall prove that if $u \geq 2$, then there are infinitely many primes $p$ such that $0 \notin f(H)$ and so $ f(H) = \Fp^{\times}.$  We begin with a standard result about irreducible polynomials in \Q.

\bl   \label{BLF:lemma:irreduc}
Let $a$ be a nonzero rational number, and let $q$ be a prime number such that $a$ is not a $q$th power.  Then the polynomial $g(x) = x^q-a$ is irreducible in $\Q[x].$
\el

\begin{proof}
Choose $\alpha \in \C$ such that $\alpha^r = a,$ and let $\zeta$ be a primitive $q$th root of unity.  Then
\[
g(x) = \prod_{i=0}^{q-1} (x-\alpha \zeta^i).
\]
If $g(x)$ factors in $\Q[x]$, then there exist polynomials $h(x), k(x) \in \Q[x]$ such that $g(x) = h(x)k(x) \in \Q[x]$ and $1 \leq \deg(h(x)) \leq q-1$. Let $r=\deg(h(x))$.  Since $q$ is prime, we have $(r,q)=1$ and there are integers $m$ and $n$ such that $rm+qn=1.$  There is also a set $I \subseteq \{0,1,2,\ldots, q-1\}$ with $|I|=r$ such that
\[
h(x) = \prod_{i\in I} (x-\alpha\zeta^i) = x^r + \cdots + (-1)^r\alpha^r\zeta^j,
\]
where $j = \sum_{i\in I}i$.  Then
\[
\beta =  (-1)^r\alpha^r\zeta^j \in \Q
\]
and
\[
\beta^m =  (-1)^{rm}\alpha^{rm}\zeta^{jm} 
=  (-1)^{rm}\alpha^{1-qn}\zeta^{jm} 
=  \frac{(-1)^{rm}\alpha\zeta^{jm}}{a^n} \in \Q.
\]
Multiplying by $(-1)^{rm}a^n$, we obtain
\[
\gamma = (-1)^{rm}a^n\beta^m = \alpha\zeta^{jm} \in \Q
\]
and so
\[
g(\gamma) = 0.
\]
Therefore, if $g(x)$ factors in $\Q[x],$ then $g(x)$ has a root $\gamma \in \Q$, 
which is impossible since the $q$th roots of $a$ are irrational.  It follows that $g(x)$ is an irreducible  polynomial.
\end{proof}

\bt       \label{BLF:theorem:mordell-uv}
Let $u$ and $v$ be relatively prime integers such that $u > |v|\geq 1$.
Let
\[
f(x,y) = ux+vy.
\]
There exist finite sets $A$ and $A'$ of integers such that
\[
|f(A)| < |s(A)|
\]
and
\[
|f(A')| < |d(A')|
\]
\et

\begin{proof}
Let $k$ be a positive integer and let $p$ be a prime such that $p \equiv 1 \pmod{k}.$  Let $H$ be the set of $k$th powers of elements of $\F_p^{\times}$.  Then $H$ is a multiplicative subgroup of $\F_p^{\times}$ of order $|H| = (p-1)/k$
and $[\F_p^{\times}:H] = k.$

Let $u$ and $v$ be integers relatively prime to $p$, and consider  the binary linear form
\[
f(x,y) = ux+vy.
\] 
By Theorem~\ref{BLF:theorem:mordell}, if $p > k^4,$ then
$\F_p^{\times} \subseteq f(H)$ and so
\[
|f(H)| \geq p-1.
\]

Let  $u = 1$ and $v = \pm 1.$  For the polynomial $d(x,y) = x-y,$ we have $d(1^k,1^k) = 0  \in d(H)$ and so $d(H) = \F_p.$
If $k$ is odd, then for the polynomial $s(x,y) = x+y$ we have $s(1^k,(-1)^k) = 0 \in s(H)$ and $s(H) = \F_p.$  In both cases,
\[
|d(H)| = |s(H)| = p.
\]

Let $u > |v| \geq 1$ and $(u,v)=(uv,p)=1$.  Let $a = -uv^{k-1}.$   If there exist $h_1 = \ell_1^k \in H$ and $h_2 = \ell_2^k \in H$ such that
\[
f(h_1,h_2) = uh_1 + vh_2 = u\ell_1^k + v\ell_2^k = 0,
\]
then 
\[
\left(\frac{u\ell_1}{\ell_2}\right)^k + u^{k-1}v = 0
\]
and the polynomial $g(x) = x^k-a$ has a root in $\F_p$, and is, therefore, reducible.   It follows that if $g(x)$ is irreducible in $\F_p[x]$, then $0 \notin f(H)$ and so $|f(H)| = p-1.$   

The rational integer $a=-u^{k-1}v$ is not a $q$th power for all sufficiently large primes $q$.  By Lemma~\ref{BLF:lemma:irreduc}, the polynomial 
$
g(x) = x^q-a
$ 
is irreducible over \Q.
Let  \mcp\ be the set of primes $p > q^4$ such that  $p \equiv 1\pmod{q}$ and $g(x)$ is irreducible in $\F_p[x]$.  The Chebotarev density theorem implies that the series $\sum_{p\in \mcp} 1/p$ diverges, and Theorem~\ref{BLF:theorem:mordell-uv} now follows directly from Theorem~\ref{BLF:theorem:rectifycong}.  This completes the proof.
\end{proof}

\def\cprime{$'$} \def\cprime{$'$} \def\cprime{$'$}
\providecommand{\bysame}{\leavevmode\hbox to3em{\hrulefill}\thinspace}
\providecommand{\MR}{\relax\ifhmode\unskip\space\fi MR }
\providecommand{\MRhref}[2]{%
  \href{http://www.ams.org/mathscinet-getitem?mr=#1}{#2}
}
\providecommand{\href}[2]{#2}

\end{document}